\documentstyle[a4]{amsart}

\def\myauthor{Ph\`ung H{{\accent"5E o}\kern-.28em\raise.2ex\hbox{\char'22}\kern-.20em} H{a\kern-.370em\raise.16ex\hbox{\char'47}\kern.1em}i}
\def\myAUTHOR{ PH\`UNG H{{\accent"5E O}\kern-.38em\raise.8ex\hbox{\char'22}\kern-.12em}  H{A\kern-.46em\raise.80ex\hbox{\char'47}\kern.18em}i}
\newcommand{\BBb}[1]{{\Bbb #1}}

\def\Mn{{M^{\ot n}}}
\def\Msn{{M^{*\ot n}}}

\def\lora{\longrightarrow}
\def\ot{\otimes} 

\def\loma{\longmapsto}
\def\Vn{{V^{\ot n}}}

\def\si{\sigma}
\newcommand{\bbas}{\begin{eqnarray*}}
\newcommand{\eeas}{\end{eqnarray*}}
\newcommand{\bbar}{\begin{array}}
\newcommand{\eear}{\end{array}}
\newcommand{\bbs}{\begin{displaymath}}
\newcommand{\ees}{\end{displaymath}}
\newcommand{\bb}{\begin{equation}}
\newcommand{\eqbb}{\begin{equation}}
\def\ee{\end{equation}}
\def\eqee{\end{equation}}
\def\eea{\end{eqnarray}}

\def\bba{\begin{eqnarray}}
\newtheorem{thm}{Theorem}
\newtheorem{lem}[thm]{Lemma}

\newtheorem{cor}[thm]{Corollary}

\newtheorem{pro}[thm]{Proposition}

\def\Im{\mbox{\rm Im}{}}

\def\Hom{\mbox{\rm Hom{}}}

\def\End{\mbox{\rm End{}}}
\newcommand{\bK}{{\Bbb K}}
\newcommand{\bZ}{{\Bbb Z}}

\def\H{{\cal H}}

\def\A{{\cal A}}

\def\P{{\cal P}}

\def\eee{\rule{.75ex}{1.5ex}\vskip1ex}

\def\proof{{\it Proof.\ }}
\newcommand{\va}{\varepsilon}

\newcommand{\rk}{{\rm rank}}
\newcommand{\tr}{{\rm tr}}
\newcommand{\Tr}{{\rm Tr}}

\def\rref#1{(\ref{#1})}

\font\Fraktur=eufm10 scaled\magstep1          % for display- and textstyle
   \newcommand{\fraktur}[1]{\mbox{\Fraktur #1}}  %
   \font\Fraktu=eufm7 scaled\magstep1            % for scriptstyle
   \newcommand{\fraktu}[1]{\mbox{\Fraktu #1}}    %
   \font\Frakt=eufm5 scaled\magstep1             % for scriptscriptstyle
  \newcommand{\frakt}[1]{\mbox{\Frakt #1}}      %
   \def\fr#1{\mathchoice{\fraktur {#1}}            % displaystyle
                        {\fraktur {#1}}            % textstyle
                        {\fraktu {#1}}             % scriptstyle
                        {\frakt {#1}}  }           % scriptscriptstyle

\newcommand{\Ss}{\fr S}

\def\db{{\mathchoice{\mbox{\rm db}}
                    {\mbox{\rm db}}
                    {\mbox{\scriptsize\rm db}}
                    {\mbox{\tiny\rm db}} }}

\def\ev{{\mathchoice{\mbox{\rm ev}}
                    {\mbox{\rm ev}}
                    {\mbox{\scriptsize\rm ev}}
                    {\mbox{\tiny\rm ev}} }}

\def\id{{\mathchoice{\mbox{\rm id}}
                    {\mbox{\rm id}}
                    {\mbox{\scriptsize\rm id}}
                    {\mbox{\tiny\rm id}} }}

\def\part{\vdash}

\newcommand{\bfLambda}{\mbox{\boldmath$\Lambda$\unboldmath}}

\newcommand{\mono}[3]{#1_{#2_1}^{#3_1}#1_{#2_2}^{#3_2}\cdots #1_{#2_n}^{#3_n}}

\def\rank{\mbox{\rm rank}}

\begin{document}

\def\mytitle{Characters of Representations of Quantum Groups of Type $A_n$}
\def\myemail{{\small\tt phung@@mpim-bonn.mpg.de}}
\def\myaddress{\small\it Max-Planck Institut f\"ur Mathematik, Gottfried-Claren-Str. 26, 53225, Bonn, Germany }
\def\myperaddress{Hanoi Institute of Mathematics P.O. Box 631, 10000 Bo Ho, Hanoi, Vietnam}
\def\mythanks{The author would like to thank the Abdus Salam International Centre for Theoretical Physics and the Max-Planck Institute for Mathematik, Bonn, for their hospitality and financial supports.}
\def\myabstract{We introduce the notion of characters of comodules over coribbon Hopf algebras. The case of quantum groups of type $A_n$ is studied. We establish a characteristic equation for the quantum matrix and a q-analogue of Harish-Chandra-Itzykson-Zuber integral}
\def\amshead{
\title{Characters of Representations of Quantum Groups of Type $A_n$}
\author{ PH\`UNG H{{\accent"5E O}\kern-.38em\raise.8ex\hbox{\char'22}\kern-.12em}  H{A\kern-.46em\raise.80ex\hbox{\char'47}\kern.18em}i}
\address{\myperaddress}
\curraddr{\myaddress}
\email{phung@@mpim-bonn.mpg.de}
\begin{abstract}\myabstract\end{abstract}
\subjclass{Primary 16W30,17B37 , Secondary 17B10}
\maketitle }
%%%%%%%%%%%%%%%%%% END  DEFINITION %%%%%%%%%%%%%%%%%%%%%%%%
\amshead
\bibliographystyle{plain}
\section*{Introduction}

Let $G$ be an algebraic group, $V$ be a $G$-module. The character $\Phi_V$ of $V$ is a function on $G$, whose value at an element $g$ of $G$ is the trace of the operator on $V$, induced by $g$. In this note we introduce the notion of characters of quantum groups. Thus the algebraic group $G$ will be replaced by a coribbon Hopf algebra $H$, $G$-module will be replaced by $H$-comodules, their characters are elements of $H$.

To the further generalization, we also introduce the notion of characters of endormorphisms of $H$-comodules, thus obtain a trace map from the endomorphism rings of $H$-comodules into $H$. Further we introduce the notion of partial characters. As an application we define the braided powers of a multiplicative matrix.

We focus ourselves on matrix quantum groups of type $A_{n-1}$. These quantum groups are defined in terms of even Hecke symmetries. They are quantum analogies ofthe general linear group $GL(n)$, whose characters are closely related to symmetric functions in $n$ variables. It turns out that the characters of a quantum group of type $A_{n-1}$ form  a ring that is isomorphic to the ring of symmetric functions in $n$ variables. The characters of endomorphisms of the tensor powers of the basic comodule induces a trace map on the Hecke algebra with values in the ring of characters. Finally, using the partial characters we define the quantum powers of a quantum matrix and prove the characteristic equation for these powers.

As an application, we give a $q$-analogue of the Harish-Chandra-Itzykson-Zuber (HCIZ) integral formula. Classically, the HCIZ formula computes the integral $I(M,N,t):=\displaystyle\int_{U(n)}\exp(t\tr(MUNU^\dag))dU$ on the unitary group $U(n)$, where $M,N$ are hermitian matrices. This integral is computed in terms of irreducible characters of $U(n)$ evaluated at $M$ and $N$ \cite{iz1}:
\bbs I(M,N,t)=\sum_{n=0}^\infty \frac{t^n}{n!}\sum_{\lambda\part n}\frac{d_\lambda}{r_\lambda}\Phi_\lambda(M)\Phi_\lambda(N).\ees
We also derive a formula for a special trace map of the Hecke algebras at primitive idempotents. Consequently, we obtaind an interesting interpretations of the characters tables of the Hecke algebras.

The paper is organized as follows. In Section \ref{sec1} we give the  definitions and basic properties of coquasitriangular and coribbon Hopf algebras. Then we define the characters of comodules over coribbon Hopf algebras as well as characters of endomorphisms and partial characters. 
In Section \ref{sec2} we first recall the definitions and representation theory of a matrix quantum group of type $A_{r-1}$. It is defined in terms of an even Hecke symmetry $R$ of rank $r$ over a vector space $V$. For more details the reader is referred to \cite{ph97,ph97b}. The Hecke symmetry induces an action of the Hecke algebra over the tensor powers of $V$, which is the centralizer of the action of the quantum group on $V$. The character map then induces a trace map on the Hecke algebras with values in the commutative algebra spanned by characters. We then define the quantum powers of the quantum matrix and prove the characteristic equation. 
In Section \ref{sec3} we prove a quantum analogue to the HCIZ formula for the integral on quantum groups of type $A_{r-1}$ in terms of the just introduced characters. At the end we also notice an interesting interpretation of the characters tables of Hecke algebras in terms of their minimal central idempotents.

 \subsubsection*{Notations}
 All bialgebras, Hopf algebras as well as their comodules are defined over a fixed  algebraically closed field $\bK$ of characteristic zero.
 For a bialgebra, the coproduct is denoted by $\Delta$, the counit is denoted by $\va$. For a Hopf algebra, the antipode is denoted by $S$.

 The Hecke algebra $\H_n=\H_{q,n}$ is spanned as a vector space over $\bK$ by a basis $\{T_w|w\in\Ss_n\}$, $\Ss_n$ is the permutation group. The element $T_{(i,i+1)}$ is denoted by $T_i$. $\H_n$ is generated, as an algebra, by $\{T_i, i=1,2,...,n-1\}$ subject to the relation $T_iT_{i+1}T_i=T_{i+1}T_iT_{i+1}$, $T_i^2=(q-1)T_i+q$. It is semisimple provided $q^n\neq 1, \forall n>1$.

A partition $\lambda$ of a positive integers $n$ is a sequence of non-increasing non-negative integers, whose sum is $n$, in notation $\lambda\part n$ of $|\lambda|=n$. $\P$ denotes the set of all partitions, $\P_n$ denotes the set of partitions of an integer $n$, $\P^r$ denotes the set of all partitions of length (i.e., the number of non-zero components) at most $r$. 

A diagram $[\lambda]$ is a matrix, whose first row contains $\lambda_1$ entries, second row contains $\lambda_2$ entries,  and so on... 
The coordinate of a node in the $i$-th row and $j$-th column is $(i,j)$.

 \section{Coribbon Hopf Algebras and Their Characters}\label{sec1}
 \subsection{Definitions}\label{sec11}
 A coquasitriangular (CQT) structure on a bialgebra $B$ is a bilinear form $r:B\otimes B\lora {\Bbb K}$ subject to the following conditions
 \begin{itemize}
 \item[i.] $a_{(1)}b_{(1)}r(a_{(2)},b_{(2)})=r(a_{(1)},b_{(1)})b_{(2)}a_{(2)}$ -- ``naturality'',
 \item[ii.] $r$ is invertible in $\Hom(B\otimes B,{\Bbb K})$,
 \item[iii.] $r(a,bc)=r(a_{(2)},b)r(a_{(1)},c)$\\
  $r(ab,c)=r(a,c_{(1)})r(b,c_{(2)})$
 -- ``multiplicativity''.\end{itemize}
 Here we use Sweedler's sigma notation for the coproduct: $\Delta(x)=x_{(1)}\otimes x_{(2)}.$

 A (Hopf) bialgebra equipped with such a structure is called coquasitriangular (CQT) (Hopf) bialgebra.

 In a CQT Hopf algebra, the square of the antipode is co-inner, that is
 \begin{equation}\label{eq0} S^2(a)=u^{-1}(a_{(1)})a_{(2)}u(a_{(3)}),\end{equation}
 where $u(a):=r(a_{(2)},S(a_{(1)})$ and $u^{-1}$ is its inverse in $\Hom(B,{\Bbb K})$, $u^{-1}(a)=r(S^2(a_{(2)}),a_{(1)})$.

 A CQT structure on a bialgebra induces a braiding \cite{js1} in the category of its right comodules in the following way. Let $(B,r)$ be a CQT bialgebra, for any right B-comodules $M,N$, define a morphism 
 \bbas \tau_{M,N}(m,n)=n_{(0)}\ot m_{(0)}r(m_{(1)},n_{(1)}):M\ot N\lora N\ot M.\eeas
 Then $\tau$ is a braiding in $B$-Comod -- the category of right $B$-comodules.

A  coribbon Hopf algebra is a CQT Hopf algebra $(H,r)$ equipped with a linear mapping $t:B\lora {\Bbb K}$ subject to the following relations:
 \begin{itemize}
 \item[i.] $t(a_{(1)})a_{(2)}=a_{(1)}t(a_{(2)})$ -- ``naturality'',
 \item[ii.] $t\circ S=t$ -- ``rigidity'',
 \item[iii.] $t$ is invertible in $\Hom(B,{\Bbb K})$,
 \item[iv.] $t(ab)=t(a_{(1)})t(b_{(1)})r(b_{(2)}a_{(2)})r(a_{(3)},b_{(3)})$ -- ``twist''.
 \item[v.] $t(1_H)=1$ -- ``normalization''.
 \end{itemize}
  $t$ has the following property:
 \begin{equation}\label{eq1.1}
 t(a_{(1)})t(a_{(2)})=u^{-1}(a_{(1)})u^{-1}(S(a_{(2)})).\end{equation}
 The proof is to write out the right-hand side of the quality $\va(a)=t(a_{(1)}S(a_{(2)}))$ and to use \rref{eq0}.

 Let $(H,r,t)$ be a coribbon Hopf algebra. Then $H$-comod, the category of finite dimensional right $H$-comodules, is a ribbon category \cite{rt,js1} with the twist given by
 \begin{equation}\label{eq01} \theta_M:M\lora M;\quad \theta_M(m)=m_{(0)}t(m_{(1)}),\end{equation}
 for any finite dimension right $H$-comodule $M$. Then $\theta$ commutes with any morphism and we have
\bbas
&\theta_{M^*}=\theta_M^*\qquad \theta_\bK=\id_\bK&\\
&\theta_{M\ot N}=(\theta_M\ot \theta_N)\tau_{N,M}\tau_{M,N}.&\eeas

\subsection{The Characters}
Now we are going to define the characters of comodules over a coribbon Hopf algebra. First we recall the definition of rank or braided dimension of comodules. The rank of a comodule $M$ over a coribbon Hopf algebra $H$ is
\bbs \rk(M):=\ev_V\circ\tau_{M,M^*}(\theta_M\ot \id_{M^*})\db_M(1_\bK)\in\bK,\ees
where $\db_M:\bK\lora M\ot M^*$, $\ev_M:M^*\ot M\lora \bK$ are morphisms that make $M^*$ the dual comodule to $M$. The isomorphism $\theta$ is inserted in the definition in order to provide that $\rk(M\ot  N)=\rk(M)\cdot \rk(N)$. It is also obvious that $\rk(M\oplus N)=\rk(M)+\rk(N)$, $\rk(M^*)=\rk(M)$. Unfortunately, the rank is generally not additive with respect to short exact sequences. For a morphism $f:M\lora M$, we define its braided trace to be
\bbs \Tr(f):=\ev_V\circ\tau_{M,M^*}(f\circ\theta_M\ot \id_{M^*})\db_M(1_\bK)\in\bK.\ees
It is easy to see that $\Tr$ is additive and tensor multiplicative, moreover, $\Tr(f\circ g)=\Tr(g\circ f)$. Notice that $\rk(M)=\Tr(\id_M)$.

The character $\Phi(M)$ is defined to be
\bbs \Phi(M):=(\ev_M\ot\id_H))(\id_{M^*}\ot \delta_M)\tau_{M,M^*}(\theta_M\ot \id_{M^*})\db_M(1_\bK)\in H.\ees
Analogously, for a morphism $f:M\lora M$, we define it character to be
\bbs \Phi(f):=(\ev_M\ot\id_H)(\id_{M^*}\ot \delta_M)\tau_{M,M^*}(f\circ\theta_M\ot \id_{M^*})\db_M(1_\bK)\in H.\ees
Notice that $\Phi(M)=\Phi(\id_M)$ and $\va(\Phi(f))=\Tr(f)$, $\va$ is the counit on $H$. Thus we have a character map $\Phi:\End^H(M)\lora H$.
\begin{lem}\label{lem1}
The character map $\Phi$ has the following properties 
\begin{itemize}\item[i.]$\Phi(f+g)=\Phi(f)+\Phi(g)$
 \item[ii.] $\Phi(f\ot h)=\Phi(f)\Phi(h)$
 \item[iii.] $\Phi(f\circ g)=\Phi(g\circ f)$.
 \end{itemize}\end{lem}
(ii) and (iii) imply that $\Phi(f)\Phi(g)=\Phi(g)\Phi(f)$.
\proof. The properties (i)and (iii) are obviously.  (ii) follows from the following identity
\bba\label{identity}\tau_{\Mn,\Msn}(\theta_{\Mn}\ot\id_\Msn)=\db_M^n\circ\db_M^{n-1}\circ\cdots\circ\db_M:\bK\lora \Mn\ot \Msn,\eea
where $\db^i_\Mn:=\id_M^{i-1}\ot\db_M\ot\id^{i-1}_{M^*}.$\eee

Thus, $\Phi$ is a trace map on $\End^H(M)$ with value in the commutative subalgebra of $H$ spanned by characters.

Assume now that $f:M\ot N\lora M\ot N$. Then we define its partial character with respect to $M$, $\Phi_M(f):N^*\ot N\lora H$, to be
\bbs \Phi_M(f):=(ev_{M\ot N}\ot \id_H)\delta_{M\ot N}(\id_{N^*}\ot(f\circ\tau_{M,M^*}(\theta_M\ot\id_{M^*})\db_M)\ot\id_N).\ees
The partial character is used particularly to define the braided powers of the multiplicative matrix associated to $M$. We mention the following property of the partial character, which will be used later. Let $f:M\lora M$ and $g:M\ot N\lora M\ot N$. Then
\bba\label{eq2}\Phi_M((f\ot\id_N)\circ g)=\Phi_M(g\circ(f\ot\id_N)).\eea

Let us now give some coordinate calculation. Fix a basis $x_1,x_2,...,x_n$ of $M$. The coaction of $H$ can then be given in terms of a multiplicative matrix $A=(a^i_j)_{i,j=1}^n$: $\delta(x_i)=x_j\ot a^j_i$, $\Delta(a^i_j)=a^i_k\ot a^k_j$, here and later on we shall adopt the convention of summing up by the indices that appear both in upper and lower places. Let $\xi_1,\xi_2,...,\xi_n$ be the dual basis on $M^*$. Thus, $\db_M(1_\bK)=x_i\ot \xi^i$, $\ev_M(\xi^i\ot x_j)=\delta^i_j.$

Assume that with respect to the bases above $\tau_{M,M}$ and $\theta_M$ have matrices $R^{ij}_{kl}$ and $T^i_j$, respectively:
\bbs \tau_{M,M}(x_i\ot x_j)=x_k\ot x_lR^{kl}_{ij},\quad \theta_M(x_i)=x_jT^j_i.\ees
Then we have $r(a^i_j,a^k_l)=R^{ki}_{jl}$, $t(a^i_j)=T^i_j$. Set $P^{ki}_{jl}=r(a^i_j,S(a^k_l)).$ Then $R^{ip}_{jq}P^{qk}_{pl}=\delta^i_l\delta^k_j$ and $P$ is the matrix of $\tau_{M,M^*}$: $\tau_{M,M^*}(x_i\ot \xi^j)=\xi^l\ot x_kP^{jk}_{il}.$ Further, set $C^i_j=u(a^i_j)$, $D^i_j=u(S(a^i_j))$. From the definition of $u$, we have $C^i_j=P^{il}_{jl}$, $D^i_j=P^{li}_{lj}$. Equation \rref{eq1.1} implies $CD=DC=T^{-2}$.

Assume that a morphism $f$ has a matrix $F$. Then
\bba\label{eq3} \Phi(f)=\tr(DTFA),\quad \Phi(f^*)=\tr(CTFA),\eea
where $\tr$ denotes the usual trace.

We define the braided power of the multiplicative matrix $A$ to be the image of $\xi^i\ot x_j$ under $\Phi_{M^{\ot n-1}}(\tau_{M^{\ot n-1},M})$. Using the third property of $\theta$, we have
\bba\label{eq4} (A^{q*n})^i_j=\Phi_{M^{\ot n-1}}(\tau_{M^{\ot n-1},M})(\xi^i\ot x_j)=(D^{\ot n-1}T^{\ot n-1})^I_J(R_{n-1}R_{n-2}\cdots R_1)^{Jj}_{Kk}A^K_Ia^k_i,\eea
where $a^I_J:=a^{i_1}_{j_1}a^{i_2}_{j_2}\cdots a^{i_n}_{j_n}.$

\subsection{Characters of Cosemisimple Coribbon Hopf Algebras}
 A Hopf algebra is called (right) cosemisimple iff the category of its (right) comodules is semisimple, that is all comodules are absolutely decomposable. In this case there exists an $H$-comodule morphism, called integral, $\displaystyle\int:H\lora {\Bbb K}$, that is, $\displaystyle\int$ satisfies the following relation: $\displaystyle\int(a)=a_{(1)}\displaystyle\int(a_{(2)}), \displaystyle\int(1)=1$. $H$ decomposes into the direct sum of its simple subcoalgebra $H=\oplus_{i\in I}H_i, H_0\cong {\Bbb K}$ and $\displaystyle\int|_{H_i}=0$ unless $i=0$. As in the case of finite group characters, we define a scalar product on the set of   characters in $H$:
 \begin{equation}\label{eq5} <\Phi(M),\Phi(N)>:=\int(\Phi(M\ot N^*))=\int(\Phi(M)\Phi(N^*)).\end{equation}
 From the preceding discussion, one sees that for a simple $H$-comodule $M$, $\displaystyle\int(\Phi(M))=0$ unless $M\cong {\Bbb K}$. Thus, we have the orthogonal relations for simple $H$-comodules
 \begin{equation}\label{eq6}
 <\Phi(M),\Phi(N)>=\left\{\begin{array}{lll}1 &\mbox{ if }& M\cong N\\
 0 & \mbox{ if } & M\not\cong N.\end{array}\right.\end{equation}
 The direct decomposition of $H$ into its simple subcoalgebras also implies that non-isomorphic simple comodules have different   characters, and the set of characters are linear independent. 
 \begin{pro}\label{dl1} Let $H$ be a co-semi-simple Hopf algebra.
 The subspace of $H$ spanned by   characters is a subalgebra of $H$ and is isomorphic to the Grothendieck ring (defined over $\BBb{K}$) of finite dimensional $H$-comodules.\end{pro}

 \section{Matrix Quantum Groups of Type $A_n$ and Their Characters}\label{sec2}
 Let $V$ be a vector space of finite dimension over ${\Bbb K}$ with a fixed basis $x_1,x_2,...,x_d$. Then an operator $R:V\ot V\lora V\ot V$ can be given by its matrix $R^{ij}_{kl}$. $R$ is called Hecke symmetry if it satisfies the Yang-Bater equation $R_1R_2R_1=R_2R_1R_2$, $R_1:=R\ot\id_V$, $R_2:=\id_V\ot R$, the Hecke equation $(R-q)(R+1)=0$ and is closed, i.e., there exists a matrix $P^{ij}_{kl}$, such that $P^{ip}_{jq}R^{qk}_{pj}=\delta^i_l\delta^k_j.$ We shall always assume that $q^n\neq 1, \forall n>1$.

 The name Hecke symmetry comes from the fact that such and operator induces a representation of the Hecke algebra $\H_n=\H_{q,n}$ on the  power $\Vn$ of $V$: $\rho_n(T_i)=R_i:=\id_V^{\ot i-1}\ot R\ot \id^{\ot n-i-1}_V.$ We denote, for convenience, $R_w:=\rho_n(T_w)$. By means of $\rho_n$ we shall some times identify an element of $\H_n$ with an endomorphism of $\End(\Vn)$.

The Hecke symmetry $R$ is called even Hecke symmetry of rank $r$ if the anti-symmetrizer operator $Y_n:=\sum_{w\in\Ss_n}(-q)^{-l(w)}R_w$ has non-zero rank for $n=r$ and vanishes for $n=r+1$.

Given a Hecke symmetry $R$, we define a bialgebra $E_R$ and a Hopf algebra $H_R$ in the following way. Let $\{z_j^i,t^i_j|i,j=1,2,...,d\}$ be a set of variables. $E_R:=\bK<\{z_j^i|i,j=1,2,...,d\}>/(R^{ij}_{mn}z^m_kz^n_l=z^i_pz^j_qR^{pq}_{kl})$, $H_R:=\bK<\{z_j^i,t^i_j|i,j=1,2,...,d\}>/(R^{ij}_{mn}z^m_kz^n_l=z^i_pz^j_qR^{pq}_{kl}, t^i_kz^k_j=z^i_kt^k_j=\delta^i_j).$ The coproduct on $E_R$ and $H_R$ is given by $\Delta(z^i_j)=z^i_k\ot z^k_j$, $\Delta(t^i_j)=t_j^k\ot t_k^i.$ The antipode on $H_R$ is given by $S(z^i_j)=t^i_j$, $S(t^i_j)=C^i_mz^m_n{C^{-1}}^n_j,$ $C^i_j:=P^{il}_{jl}$.

$E_R$ and $H_R$ are coquasitriangular bialgebra with the CQT structure given by $r(z^i_j,z^k_l)=R^{ki}_{jl}$. $H_R$ is a coribbon Hopf algebra with the coribbon structure $t(z^i_j)=q^{(r+1)/2}\delta^i_j.$

$E_R$ and $H_R$ coact on $V$ by the coaction $\delta(x_i)=x_j\ot z^j_i.$ Hence $\Vn$ is a comodule of $E_R$ and $H_R$, its dual is comodule of $H_R$. The coaction of $E_R$ and the action of $\H_n$ on $\Vn$ are centralizers of each other in $\End_\bK(\Vn)$. Thus, $\rho_n(\H_n)\cong\End^{E_R}(\Vn)$ and simple $E_R$ comodules have form $\Im(\rho_n(E_\lambda))$, where $E_\lambda$ is a primitive idempotent of $\H_n$. The natrural map $E_R\lora H_R$ is injective, hence simple $E_R$-comodule are simple over $H_R$, too.

If $R$ is an even Hecke symmetry of rank $r$ then $H_R$ is cosemisimple, its simple comodules are parameterized by $\bZ$-partitions, that is, the sequences of non-increasing integers of length $r$: $\lambda=(\lambda_1,\lambda_2,...,\lambda_r), \lambda_i\in\bZ, \lambda_i\geq \lambda_{i+1}.$ Let $M_\lambda$ denote the simple comodule corresponding to $\lambda$. In particular, $V\cong M_{(1)}$. The coefficients $ c_{\mu\nu}^\lambda $ in the decomposition
 \begin{equation}\label{eq61}M_\mu\ot M_\nu\cong \bigoplus_\lambda c_{\mu\nu}^\lambda M_\lambda\end{equation}
 are the same as for those of irreducible rational GL$(r)$-representations, i.e., the Littlewood-Richardson coefficients. Further, let $-\lambda:=(-\lambda_r,-\lambda_{r-1},...,-\lambda_1),$ $\lambda^+:=(\lambda_1-\lambda_r,\lambda_2-\lambda_r,...,0).$ Then
\bbs (M_\lambda)^*\cong M_{-\lambda},\quad M_\lambda\cong M_{\lambda^+} \ot M{(1^r)}^{\ot \lambda_r}.\ees
Notice that $M_{(1^n)}$ is isomorphic to the image of $Y_n$ in $\Vn$. It is shown that $Y_r$ has rank 1 over $\bK$ hence induces a group-like element $D$ in $E_R$, called quantum determinant. 

The interested reader is referred to \cite{ph97,ph97b} for more details.

 \subsection{Characters of Comodules over $H_R$}\label{sec31} 
Now, using the construction of Section \ref{sec1} we can define characters of simple $H_R$-comodules. Notice that the matrix $C$ introduced above is precisely $C^i_j=u(z^i_j)$. We also set $D^i_j=u(t^i_j)=P^{li}_{lj}$. Then, according to \rref{eq1.1}, $CD=DC=q^{-(r+1)}.$ 

Let $\lambda\part n$. Then the simple comodule $M_\lambda$ is isomorphic to $\Im(\rho_n(E_\lambda))$, where $E_\lambda$ is a primitive idempotent corresponding to $\lambda$. Let us use the same character ($E_\lambda$) to denote the endomorphism on $\Vn$ induced by $E_\lambda$. Then we have
\begin{lem}
\bbas\Phi(M_\lambda)=q^{n(r+1)/2}\tr(D^{\ot n}E_\lambda Z^{\ot n})&& 
\Phi(M_\lambda^*)=q^{n(r+1)/2}\tr(C^{\ot n}E_\lambda \overline{T^{\ot n}})\\
\Phi(F)=q^{n(r+1)/2}\tr(D^{\ot n}F Z^{\ot n})&& 
\Phi(F^*)=q^{n(r+1)/2}\tr(C^{\ot n}F \overline{T^{\ot n}}),\eeas
where $F$ is the matrix of a morphism $f:\Vn\lora \Vn,$ $\overline{T^{\ot n}}^I_J={T^{\ot n}}^{I'}_{J'}.$
\end{lem}
\proof The above formula follow immediately from the identity \rref{identity}
Recall that $\Phi$ is a trace map from $\End^{E_R}(\Vn)$ into $H_R$. Combining $\Phi$ with $\rho_n$ we get a trace map $\H_n\lora H_R$. Using Lemma \ref{lem1} we can represent the character of a comodule morphism on $\Vn$ as a linear combination of $\Phi(T_{c_k})$, where $c_k:=(1,2,...,)$ a $k$-cycle in $\Ss_n$ (see \cite{ph97}). Let $P_n:=\Phi(T_{c_n}),$ $H_n:=\Phi(M_{(n)})$, $E_n:=\Phi(M_{(1^n)})$.
 \begin{lem}\label{lem2} We have the following relations
 \bbas [n]_qH_n&=&H_{n-1}P_1+H_{n-2}P_2+\cdots +H_0P_n,\\
{[n]}_{1/q}E_n&=&E_{n-1}P_1-q^{-1}E_{n-2}P_2+\cdots +(-q)^{n-1}E_0P_n ,\eeas
 where $ [n]_q:=(q^n-1)/(q-1)$.\end{lem}
 The second relation was established in \cite{ps96} for  slightly differently defined $P_n$'s.

 \proof Let $X_n$ and $Y_n$ be the quantum symmetrizer and anti-symmetrizer operators on $V^{\ot n}$. They can be defined by induction:
 \bbas {[n]}_qX_n&=&(1+R_{n-1}+\cdots +R_{1}R_2\cdots R_{n-1})X_{n-1},\\
{[n]}_{1/q} Y_n&=&(1-q^{-1}R_{n-1}+\cdots +(-q)^{n-1}R_{1}T_2\cdots R_{n-1})Y_{n-1}.\eeas
  $X_n$ and $Y_n$ are projectors on $V^{\ot n}$ and their images are isomorphic to $M_{(n)}$ and $M_{(1^n)}$ respectively. Hence $\Phi(X_n)=S_{(n)}=H_n$ and $\Phi(Y_n)=S_{(1^n)}=E_n$. We have $X_nR_i=R_iX_n=qX_n, Y_nR_i=R_iY_n=-Y_n$. Therefore
\bbas [n]_qH_n&=& \Phi((1+R_{n-1}+\cdots+R_1R_2\cdots R_{n-1})X_{n-1})\\
&=& \Phi(X_{n-1}(1+R_{n-1}+\cdots+R_1R_2\cdots R_{n-1}))\\
&=& H_{n-1}P_1+\Phi((1+R_{n-2}+\cdots+R_1R_2\cdots R_{n-2})X_{n-2}R_{n-1})\\
&=&H_{n-1}P_1+\Phi(X_{n-2}(1+R_{n-2}+\cdots+R_1R_2\cdots R_{n-2})R_{n-1})\\
&=&\cdots\\
&=& H_{n-1}P_1+H_{n-2}P_2+\cdots+P_n.\eeas
The second equation is proved analogously.\eee
\begin{cor}The character map $\Phi$ is a trace map on $\H_n$ with value in $\Lambda_R$, the subalgebra of $E_R$ spanned by $\{ S_\lambda:=\Phi(M_\lambda)|\lambda\in\P^r\}$.\end{cor}

On the other hand, according to Equation \rref{eq6} and Lemma \ref{lem1}, we have
\begin{equation}\label{detform} S_\lambda=\mbox{det}\left|(H_{\lambda_i-i+j})_{1\leq i,j\leq r}\right|,\lambda\part n\quad S_\lambda=S_{\lambda^+}\cdot D^{\lambda_r}.\end{equation} 
Therefore we can define an algebra isomorphism $\Lambda_E\lora \bfLambda_r$, $H_k\loma h_k$, where $\bfLambda_r$ is the algebra of symmetric functions in $r$ variables and $h_k$ are $k$-th complete symmetric function. This isomorphism is in fact an isometry with respect to the scalar product defined in \rref{eq3}, under which $S_\lambda$ is mapped to $s_\lambda$ -- the Schur functions and, in particular, $E_k:=S_{(1^k)}$ is mapped to the elementary symmetric function $e_k$. Thus $S_\lambda$ can be consider as quantum symmetric functions on non-commuting variables. Notice that the elements $P_k$ are mapped to the symmetric functions $q_k$, introduced by A. Ram in \cite{ram91}.

 Let $K$ be  an algebra over $\bK$. A $K$-point $\A$ of $E_R$ is an algebra homomorphism $\A:E_R\lora K$. Let $A=\A(Z)$. Then a satisfies the relation $RA_1A_2=A_1A_2R$. Conversely, each matrix of elements from $K$, satisfying this relation give a $K$-point of $E_R$, it is then called a quantum matrix with coefficients from $K$. In this case we can define the elements $S_\lambda(A)$ of $K$, replacing $Z$ by $A$ in the definition of $S_\lambda$. $S_\lambda(A)$ is called the value of $S_\lambda$ at the point $\A$. Analogously, let $E_R^{\rm op}$ be the subbialgebra of $H_R$, generated by $t$'s. Let $N$ be a point of $E_R^{\rm op}$, we define $S_{-\lambda}(N)$ to be $S_{-\lambda}$ computed at $T=N$.

Let us consider an example. Assume that $R$ is the Drinfel'd-Jimbo matrix of type $A_{r-1}$, for $1\leq i,j\leq r$, $p^2=q$,
\bbs R^{kl}_{ij}:=\left\{
\bbar{lll}
 p^2 &\mbox{ if }& i=j=k=l\\
 p^2-1 &\mbox{ if }&k=i<j=l\\
 p&\mbox{ if }&k=j\neq i=l\\
 0&\multicolumn{2}{l}{\mbox{ otherwise }}
\eear\right.\ees
Then $R$ is an even Hecke symmetry of rank $r$. The associated Hopf algebra $H_R$ is called (the function algebra on) quantum general linear group $GL_q(r)$.
The elementary symmetric functions are
\bbs E_k=p^{k(r+1)}\sum_{i_1<i_2\cdots<i_k}p^{2(i_1+i_2+\cdots i_k)}\sum_{\si\in\Ss_k}(-p)^{-l(\si)}z_{i_1}^{i_{1\si}}z_{i_2}^{i_{2\si}}\cdots z_{i_k}^{i_{k\si}}.\ees

We see that $R_{ij}^{kl}=0$ unless $(i,j)=(k,l)$ or $(i,j)=(l,k)$. Therefore any diagonal matrix with commuting entries in the diagonal is a point of GL$_q(r)$. Let $A=(a_1,a_2,\ldots, a_r)$ be such a point. Then
\bbs S_\lambda(A)=s_\lambda(qa_1,q^2a_2,\ldots, q^ra_r),\ees
where $s_\lambda$ are the usual Schur functions.

\subsection{The Characteristic Equation}

Using the construction of Section \ref{sec1} we can also define the quantum powers of the matrix $Z$.
\bbs (Z^{q*n})^i_j:=q^{n(r+1)/2}{D^{\ot n}}^M_N(R_nR_{n-1}\cdots R_1)^{Nj}_{Pq}Z^P_Mz^q_i
\ees

\begin{thm}(The characteristic equation for quantum matrix) The quantum matrix $Z$ satisfies the following equation
\begin{equation}\label{eq9}
 Z^{q*r}-E_1Z^{q*(r-1)}+\cdots+(-1)^rE_r=0.\end{equation}\end{thm}
\proof 
  If $R$ is an even Hecke symmetry of rank $r$, $Y_{r+1}=0$. Taking the partial trace of $Y_{r+1}$ with respect to $V^{\ot r}$, we have
 \bbas 0&=&\Phi_{V^{\ot r}}\left((1-q^{-1}R_{r}+\cdots +(-q)^{-r}R_{1}R_2\cdots R_{r})Y_{r}\right)\eeas

Using the equation \rref{eq2}, we have
\bbas \lefteqn{q^{n(r+1)/2}\Phi_{V^{\ot r}}((1-q^{-1}R_r+\cdots+(-q)^{-r}R_1R_2\cdots R_r)Y_r)}\\
&=& \Phi_{V^{\ot r}}(Y_r(1-q^{-1}R_r+\cdots+(-q)^{-r}R_1R_2\cdots R_r))\\
&=&\Phi_{V^{\ot r}}(Y_r)-[r]_{1/q}\Phi_{V^{\ot r}}(Y_rR_r)\\
&=&q^{-n(r+1)/2}E_r\cdot Z-\Phi_{V^{\ot r}}((1-q^{-1}R_{r-1}+\cdots+(-q)^{1-r}R_1R_2\cdots R_{r-1})Y_{r-1}R_r)\\
&=&\cdots\\
&=&q^{-n(r+1)/2}(E_r\cdot Z-E_{r-1}\cdot Z^{q*2}+\cdots+(-1)^rZ^{q*(r+1)}).\eeas
Thus, \bba\label{eq10}E_r\cdot Z-E_{r-1}\cdot Z^{q*2}+\cdots+(-1)^rZ^{q*(r+1)}=0.\eea

 We have the following equalities
 \begin{eqnarray*}
 (Z^{q*n})_i^jt^i_mP^{mk}_{jn}{D^{-1}}^n_l&=&q^{(r+1)/2}(Z^{q*(n-1)})_l^k,\\
  z^j_it^i_mP^{mk}_{jn}{D^{-1}}^n_l&=&q^{(r+1)/2}\delta^k_l.
 \end{eqnarray*}
Multiplying the left-hand side of \rref{eq10} with $t^i_mP^{mk}_{jn}{D^{-1}}^n_l$, we obtain
\bbs
 Z^{q*r}-E_1Z^{q*(r-1)}+\cdots+(-1)^rE_r=0.\ees
 {\it A priori} this equation holds in $H_R$. However, since the canonical map $E_R\lora H_R$ is injective, this equation also holds in $E_R$.\eee

 \section{The $q$-analogue of Harish-Chandra-Itzykson-Zuber integral}\label{sec3}

The HCIZ integral on $U(n)$, $\displaystyle\int e^{t\tr(MUNU^\dag)}dU,$ for any Hermitian matrices $M$ and $N$, is given by
\bbs \int e^{t\tr(MUNU^\dag)}dU=\sum_{n=1}^\infty\frac{t^n}{n!} \sum_{\lambda\part n}\frac{d_\lambda}{r_\lambda} \Phi_\lambda(M)\Phi_\lambda(N).\ees
Here $d_\lambda$ is the dimension of the irreducible representation $S_\lambda$ of the symmetric group $\Ss_n$  and $r_\lambda$ is the dimension of the irreducible  representation$ M_\lambda$ of $U(n)$, $\Phi_\lambda$ is the character of $M_\lambda$, thus $\Phi_\lambda(M)$ is the symmetric function $s_\lambda$ computed at the eigenvalues of $M$. Equivalently, we have
\bbs \int\tr(MUNU^\dag)^ndU= \sum_{\lambda\part n}\frac{d_\lambda}{r_\lambda} \Phi_\lambda(M)\Phi_\lambda(N).\ees
In this section we want to give an analogue of this formula for the integral on $H_R$.

 A formula for the normalized integral on $H_R$ is given in \cite{ph97b}. Since $H_R$ is coquasitriangular, its elements can be represented as linear combinations of elements of the form $\mono{z}{i}{j}t_{k_1}^{l_1}t_{k_2}^{l_2}\cdots t_{k_m}^{l_m}$. For simplicity we denote $Z_I^J:=\mono{z}{i}{j}$ and so on. The integral on $Z_I^JT_K^L$ vanishes unless $l(I)=l(K)$, where $l(I)$ means the length of $I$. If $l(I)=l(K)=n$, then
 \begin{equation}\label{ieq1}
\int Z_I^JT_K^L=\sum_{w\in\Ss_n}( P_nC^{\ot n}R_{w^{-1}})^{L'}_I {T_w}^J_{K'}\end{equation}
 where, for $K=(k_1,k_2,\ldots,k_n)$, $K':=(k_n,k_{n-1},\ldots, k_1)$.
\bbs P_n:=\prod_{k=1}^n(L_k-[-r]_q)^{-1}, \quad L_k:=\sum_{i=1}^{n-1}q^{-i}R_{(n-i,i)};\quad L_1:=0.\ees
 $L_k$ is the Murphy operator introduced in \cite{dj2}. Let $F_\lambda$ be the minimal central idempotent in $\H_n$, corresponding to $\lambda\part n$. Then 
\bbs P_n=\sum_{\lambda\part n}q^{rn}\prod_{x\in[\lambda]}[c_\lambda(x)+r]_q^{-1}F_\lambda,\ees
$c_\lambda(x)$ is the content of the node $x$ in the diagram $[\lambda]$ associated to $\lambda$, if $x$ is in the $i^{\rm th}$ row and  $j^{\rm th}$ column then $c_\lambda(x)=j-i.$ Let us denote $p_\lambda:=q^{rn}\prod_{x\in[\lambda]}[c_\lambda(x)+r]_q^{-1}$.

$\H_n$ possesses a non-degenerate symmetric associative bilinear form, given by $<T_u,T_w>:=\delta_u^{v^{-1}} q^{l(u)}$. Then $\chi(h)=<h,1>$ is a trace map on $\H_n$. By definition, $\chi(h)$ is the coefficient of $1$ in the presentation of $h$ as the linear combination of $T_w$. We wish to find $\chi(E_\lambda)$.

Let $\H_n=\bigoplus_{\lambda\in\P_n}\A_\lambda$ be the decomposition of $\H_n$ into the direct sum of its minimal two-sided ideals. $\A_\lambda$ is itself a matrix ring of degree $d_\lambda$, with the unit $F_\lambda$. Let $E_\lambda^{ij}$ be a $\Bbb K$ basis of $\A_\lambda$ such that $E_\lambda^{ij}E_\lambda^{kl}=\delta_k^jE_\lambda^{il}$. From standard argument (cf. \cite{g-h-j})
for any trace map $h$ on $\H_n$, $h(E_\lambda^{ij}=0$ unless $i=j$ and $h(E_\lambda^{ii}=d^{-1}_\lambda \chi(F_\lambda)=:k_\lambda$ does depend only on $\lambda$.Then
$ <E_\lambda^{ij},E_\mu^{kl}>=\delta^\mu_\lambda\delta^k_j\delta_i^lk_\lambda.$
 Therefore $\{\frac{1}{k_\lambda}E^{ji}_\lambda|\lambda\in\P_n,i,j=1,2,\ldots,d_\lambda\}$ is the dual basis to $\{E^{ij}_\lambda|\lambda\in\P_n,i,j=1,2,\ldots,d_\lambda\}$ with respect to the bilinear form defined above. On the other hand,  with respect to this bilinear form, $\{ q^{-l(w)}T_{w^{-1}}|w\in\Ss_n\}$ is the dual basis to $\{T_{w}|w\in\Ss_n\}$.
Thus, we have
\bba\label{ieq5}\sum_{w\in\Ss_n}q^{-l(w)}T_{w^{-1}}\ot T_w=\sum_{\lambda,i,j}
\frac{1}{k_\lambda}E^{ij}_\lambda \ot E^{ji}_\lambda.\eea

Now assume that $M$ is a point of $E_R$ and $N$ is a point of $E_R^{\rm op}$.
According to \rref{ieq1}, \rref{ieq5}, we have
\bbas \lefteqn{\int q^{n(r+1)/2}\tr(D^{\ot n}M^{\ot n}Z^{\ot n} \overline{N^{\ot n}}\overline{T^{\ot n}})}&&\\
&=&\sum_{w\in\Ss_n}q^{-l(w)}q^{n(r+1)/2}\tr(P_nC^{\ot n}R_{w^{-1}}\overline{N^{\ot n}})\cdot \tr(R_wD^{\ot n}{M^{\ot n}})\\
&=& \sum_{1\leq i,j\leq d_\lambda\atop{\lambda\part n}}k_\lambda^{-1}
q^{n(r+1)/2}\tr(P_nE_\lambda^{ij}C^{\ot n}\overline{N^{\ot n}})\cdot \tr(E_\lambda^{ij}D^{\ot n}{M^{\ot n}})\\
&=&\sum_{1\leq i,j\leq d_\lambda\atop{\lambda\part n}}k_\lambda^{-1}p_\lambda
q^{-n(r+1)/2}\Phi(E_\lambda^{ij*})(N)\cdot \Phi(E_\lambda^{ji})(M)\\
&=&\sum_{\lambda\part n}q^{-n(r+1)/2}d_\lambda k_\lambda^{-1}p_\lambda
S_\lambda(M)S_{-\lambda}(N).\eeas
Here ${\overline{N^{\ot n}}}^I_J:={M^{\ot n}}^{I'}_{J'}.$ In the last equation we use the fact that $\Phi(E_\lambda^{ji})=0$ unless $i=j$. Thus, to compute the integral above, it remains to find $k_\lambda$.

Since $\H_n$ is semi-simple, $\Phi$ can be uniquely represented as a combination of irreducible characters of $\H_n$ with coefficients from $\bfLambda_E$:
  $\Phi=\sum_\lambda K_\lambda\chi^\lambda$,
 where $\chi^\lambda$ is the irreducible character of $\H_n$, that corresponds to $\lambda\in\P_n$. Let $E_\lambda$ be a primitive idempotent corresponding to $\lambda$, then $\chi^\lambda(E_\mu)=\delta^\mu_\lambda$ and $\Phi(E_\lambda)=S_\lambda$, hence $K_\lambda=S_\lambda$. Thus, we have
$ \Phi=\sum_\lambda\chi^\lambda S_\lambda.$
  The orthogonal relations in \rref{eq3} imply
 \bbs \chi^\lambda(W)=(\Phi(W),\Phi(E_\lambda))=\int(\Phi(W)\Phi(E^*_\lambda)).\ees
 Recall that  (cf. \rref{eq3})
 $ \Phi(W)=q^{n(r+1)/2}(D^{\ot n}WZ^{\ot n}),$ $ \Phi(W^*)=q^{n(r+1)/2}(C^{\ot n}WZ^{\ot n}).$
 Therefore using \rref{ieq1}, \rref{ieq5}, we have
 \bbas
 \chi^\lambda(W)&=&\sum_{w\in \Ss_n}q^{-l(w)}\tr(C^{\ot n}E^{ii}_\lambda T_wWP_n T_{w^{-1}})\\
&=&k_\lambda^{-1}p_\lambda\tr(C^{\ot n}WF_\lambda),\eeas
here we use the fact that $CD=q^{-(r+1)}$.
In particular, $1=\chi^\lambda(E^{ii}_\lambda)=k^{-1}_\lambda p_\lambda\tr(C^{\ot n}E_\lambda^{ii}).$

Notice that $q^{r(n+1)/2}\tr(C^{\ot n}W)=\va(\Phi(W^*))$ is the braided trace of the morphism $W$ in $\End(\Vn)$. In particular, $q^{r(n+1)/2}\tr(C^{\ot n}E_\lambda^{ii})$ is the rank of the simple comodule $M_\lambda$, it was computed in \cite{ph97b}:
\bbs \rk(M_\lambda)=q^{-n(r-1)/2+n(\lambda)}\prod_{x\in[\lambda]}\frac{[c_\lambda(x)+r]_q}{[h_\lambda(x)]_q}.\ees
Therefore 
\bbs k_\lambda=q^{n(\lambda)}\prod_{x\in[\lambda]}{[h(x)]_q}^{-1}.\ees

Thus we have proved
\begin{pro} The value of the character $\chi$ on a primitive idempotent $E_\lambda$ is equal to $ q^{n(\lambda)}\prod_{x\in[\lambda]}{[h(x)]_q}^{-1}$, or equivalently $\chi=\sum_\lambda q^{n(\lambda)}\prod_{x\in[\lambda]}{[h(x)]_q}^{-1}\chi^\lambda.$\end{pro}
The value of $k_\lambda$ gives us immediately the quantum analogue of HCIZ integral.
\begin{thm} The  quantum  HCIZ integral on quantum group of type $A_{r-1}$ is given by
\bbs \int q^{n(r+1)/2} \tr(D^{\ot n}M^{\ot n}Z^{\ot n} \overline{N^{\ot n}}\overline{T^{\ot n}})=\sum_{\lambda\part n}\frac{d_\lambda}{\rk(M_\lambda)}S_\lambda(M) S_{-\lambda}(N).\ees\end{thm}

\subsubsection*{Example}  Assume that $R$ is the Drinfel'd-Jimbo matrix of type $A_{r-1}$ and $M=(\mu_1,\mu_2,\ldots,\mu_r)$, $N=(\nu_1,\nu_2,\ldots,\nu_r)$ be diagonal matrices. Then we have
\bbs\int(q^{n(r+1)/2} q^{|I|} \mu^IZ^I_J\nu^JT^J_I)=\sum_{\lambda\part n}\frac{d_\lambda}{\rk(M_\lambda)}s_\lambda(q^i\mu_i)s_\lambda(q^i\nu_i),\ees
where  $|I|:=i_1+i_2+\cdots+i_m$, $\mu^I:=\mu^{i_1}\mu^{i_2}\cdots \mu^{i_m}$.

\vskip2ex
Notice that the value of $k_\lambda$ allows us to give a formula for computing irreducible characters of the Hecke algebra $\H_n$. In fact, since $\chi$ is faithful, an irreducible character $\chi^\lambda$ can be computed from $\chi$ in the following sense
$ \chi^\lambda(W)=\chi(WF_\lambda)\chi(F_\lambda)^{-1},$
where $F_\lambda$ is the minimal central idempotent corresponding to $\lambda$. If $\rank(R)=n$ then $\H_n\cong \End^{E_R}(\Vn)$. Therefore we have
\begin{pro} The minimal central idempotents $F_\lambda$ contain all information about the character table of the Hecke algebra. Let $t_\lambda(w)$ be the coefficient of $T_w$ in $F_\lambda$. Then
\bbs \chi^\lambda(T_w)=\frac{q^{l(w)}q^{-n(\lambda)}}{n!}\prod_{x\in \lambda}(h(x)[h(x)]_q) t_\lambda(w^{-1}).\ees
\end{pro}
Indeed, $\chi(F_\lambda)=d_\lambda\chi(E_\lambda) = q^{n(\lambda)} d_\lambda \prod_{x\in[\lambda]} [h(x)]_q^{-1}$. $d_\lambda=n!\prod_{x\in[\lambda]}h(x)^{-1}.$ 

\vskip4ex
\begin{center} \bf Acknowledgment\end{center}

\mythanks

 \end{document}